\documentclass{article}

\usepackage{subfigure}
\usepackage{mathrsfs}
\usepackage{eucal}
\usepackage{amsmath}
\usepackage{amssymb}
\usepackage{amsfonts}
\usepackage{graphics,graphicx}
\newcommand{\Eref}[1]{Equation (\ref{#1})}
\newcommand{\fref}[1]{Figure (\ref{#1})}
\newcommand{\frefs}[1]{Figures (\ref{#1})}

\newcommand{\capn}{\mathbf{N}}

\newcommand{\xx}{\mathbf{x}}

\newtheorem{def1}{Remark}
\begin{document}
%\NME{1}{6}{00}{28}{09}

%\runningheads{Bordas \textit{et al.}} {Approximation in the SFEM}

\begin{center}
\large{\bf On the approximation in the smoothed finite element method (SFEM)}

St\'ephane~P.A.~Bordas$^{1,\ast}$,Sundararajan~Natarajan$^2$ \\

$^{1,\ast}$Lecturer, Department of Civil Engineering, University of Glasgow,
G12 8LT, Scotland, U.K. stephane.bordas@alumni.northwestern.edu. \\
$^2$PhD Research Student, Department of Civil Engineering,
University of Glasgow, G12 8LT, Scotland, U.K.
\end{center}

\begin{abstract}
This letter aims at resolving the issues raised in the recent short communication \cite{Zhang2008} and
answered by \cite{Liu2009} by proposing a systematic approximation scheme based on non-mapped shape
functions, which both allows to fully exploit the unique advantages of the smoothed finite element method
(SFEM) \cite{Liu2007,Hung2009,liu:sfem2,hbb_shell,Nguyen2007,hbb0,hbb_plate} \emph{and} resolve
the existence, linearity and positivity deficiencies pointed out in \cite{Zhang2008}.

We show that Wachspress interpolants~\cite{Wachspress1975} computed
in the physical coordinate system are very well suited to the SFEM,
especially when elements are heavily distorted (obtuse interior
angles). The proposed approximation leads to results which are
almost identical to those of the SFEM initially proposed in
\cite{Liu2007}.

These results that the proposed approximation scheme forms a strong
and rigorous basis for construction of smoothed finite element
methods.
\end{abstract}

{\bf keywords:} Smoothed Finite Element Method, boundary integration,
Wachspress Interpolants, strain smoothing, rational basis finite
elements, SFEM, isoparametric

\section{INTRODUCTION}

The smoothed finite element method (SFEM) was first proposed in \cite{Liu2007}. This new numerical method, based on gradient (strain) smoothing, is rooted in meshfree stabilized conforming nodal integration \cite{Chen2000} and was shown to provide a suite of finite elements with a range of interesting properties. Those properties depend on the number of smoothing cells employed within each finite element (see \cite{Bordas2008} for a review of recent developments and properties) and include:

\begin{itemize}
\item improved dual accuracy and superconvergence;
\item relative insensitivity to volumetric locking;
\item relative insensitivity to mesh distortion;
\item softer than the FEM.
\end{itemize}

A rigorous theoretical framework was provided in \cite{liu:sfem2,hbb0} and the method was extended to plates \cite{hbb_plate}, shells \cite{hbb_shell} and coupled with the extended finite element method \cite{Bordas2008}.
%, [include all the other ones with adequate citations].

The essential feature of the SFEM is that no isoparametric mapping is required, which implies that the approximation can be defined in the physical space directly, thereby providing freedom in the selection of the element geometry.

\begin{figure} % Figure with averaged shape functions
\centering
\scalebox{1.0}{\input{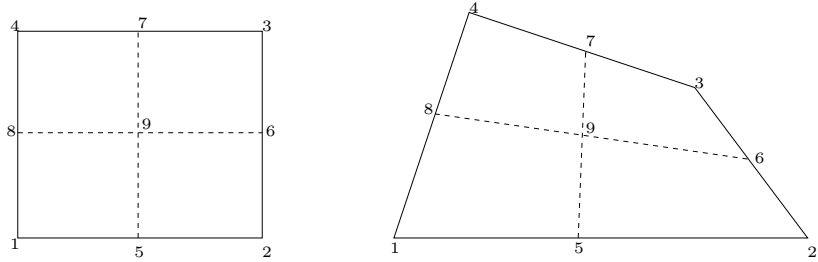}}
\caption{A four-node element divided into four smoothing cells.}
\label{fig:ShapeFunctions}
\end{figure}

In the initial paper \cite{Liu2007} (Eq. (22) and reproduced here for simplicity \Eref{eqn:liueq22}), non-mapped Lagrange shape functions are proposed as a possibility to calculate the shape functions at an arbitrary point within a smoothed finite element. It is then stated in the same paper (p863 last paragraph) that ``unless state otherwise, we still use the averaged shape functions \cite{Liu2007} for convenience.'' These shape functions are recalled in Table~\ref{table:ShapeFunctions} and \fref{fig:ShapeFunctions}, for ease of reading.

\begin{equation}
\renewcommand{\arraystretch}{1.5}
N_e(\xx_e) = \left[ \begin{array}{cccc}1 & x_e & y_e & x_e y_e \end{array} \right] \left[ \begin{array}{cccc}1 & x_1 & y_1 & x_1 y_1 \\ 1 & x_2 & y_2 & x_2 y_2 \\ 1 & x_3 & y_3 & x_3 y_3 \\ 1 & x_4 & y_4 & x_4 y_4 \end{array} \right]^{-1}
\label{eqn:liueq22}
\end{equation}

\begin{table}
\renewcommand{\arraystretch}{1}
\caption{Shape function value at different sites within an element (\fref{fig:ShapeFunctions})}
\centering
\begin{tabular}{llllll}
\hline
Site & Node 1 & Node 2 & Node 3 & Node 4 & Description \\
\cline{1-6}
\hline
1 & 1.0 & 0.0 & 0.0 & 0.0 & Field node \\
2 & 0.0 & 1.0 & 0.0 & 0.0 & Field node \\
3 & 0.0 & 0.0 & 1.0 & 0.0 & Field node \\
4 & 0.0 & 0.0 & 0.0 & 1.0 & Field node \\
5 & 0.5 & 0.5 & 0.0 & 0.0 & Side midpoint \\
6 & 0.0 & 0.5 & 0.5 & 0.0 & Side midpoint \\
7 & 0.0 & 0.0 & 0.5 & 0.5 & Side midpoint \\
8 & 0.5 & 0.0 & 0.0 & 0.5 & Side midpoint \\
9 & 0.25 & 0.25 & 0.25 & 0.25 & Intersection of two bimedians\\
\hline
\end{tabular}
\label{table:ShapeFunctions}
\end{table}

In fact, in our work on the SFEM \cite{Hung2009,hbb_shell,Nguyen2007,hbb0,hbb_plate}, and, to our knowledge, in all other work published to date \cite{Liu2007,Hung2009,liu:sfem2,hbb_shell}, these ``averaged shape functions'' have been used, with good results.

Yet, \cite{Zhang2008} provides a critique of the SFEM stating that the approximation provided by \Eref{eqn:liueq22} are inadequate because:

\begin{itemize}
\item they do not always exist (as described in the 1975 book \cite{Wachspress1975});
\item they may not be positive everywhere in the element;
\item they may not be linear everywhere in the element.
\end{itemize}

Because of this \cite{Zhang2008} disqualifies the current version of the SFEM and discredits the existing results of \cite{Liu2007,Hung2009,liu:sfem2,hbb_shell,Nguyen2007,hbb0,hbb_plate}, despite the fact (also noted in \cite{Liu2009}) that those non-mapped Lagrange shape functions of \Eref{eqn:liueq22} were in general not used in the aforementioned papers.

In this contribution, we show that it is possible to resolve the three issues mentioned by \cite{Zhang2008} about the Lagrange non-mapped shape functions while retaining the advantageous features of the smoothed finite element method, in particular its ability to deal with extremely distorted meshes.

%==================================================================
%=============== End of Introduction
%==================================================================
\section{WACHSPRESS INTERPOLANTS}
\label{section:wachspress}
Wachspress~\cite{Wachspress1975} presented a rigorous formulation for generating shape functions on arbitrary polygons, based on projective geometry. The Wachspress shape functions are unconventional compared to the polynomials used in the finite element literature, as they are in general rational functions, i.e, the ratio of two polynomials~\cite{Dasgupta2003}. These shape functions have the following essential features of interpolants (for arbitrary $n$-sided polygons):
\begin{itemize}

\item The $N_i$ satisfy the partition of unity and Kronecker Delta properties;
\item Shape function, $N_i$ is linear on sides adjacent to node $i$.
%\item Shape function, $N_i(x,y)$ vanishes on sides opposite node $i$ and at all nodes $j$ for which $i \ne j$.
\item The $N_i$ are linear complete.
%\item There is a node at each vertex and on each conic side. For each node there is an associated wedge within each polycon containing the node; % SPAB Can we make this clearer, this is a bit obscure
\end{itemize}

These properties make Wachspress interpolants  effective to build
approximations on arbitrary $n$-gons, quadrilaterals in particular.
Consider the quadrilateral (Q4), on which Wachspress interpolants
will be formulated, shown in~\fref{fig:quad}. Let
$l_{1},l_{2},l_{3}$ and $l_{4}$ be the equations of the lines
corresponding to each of the four sides of the quadrilateral, $(1-2,
2-3, 3-4 ~\&~ 4-1)$ respectively, and written in parametric form as:

\begin{figure} % Figure with averaged shape functions
\centering
\scalebox{0.6}{\input{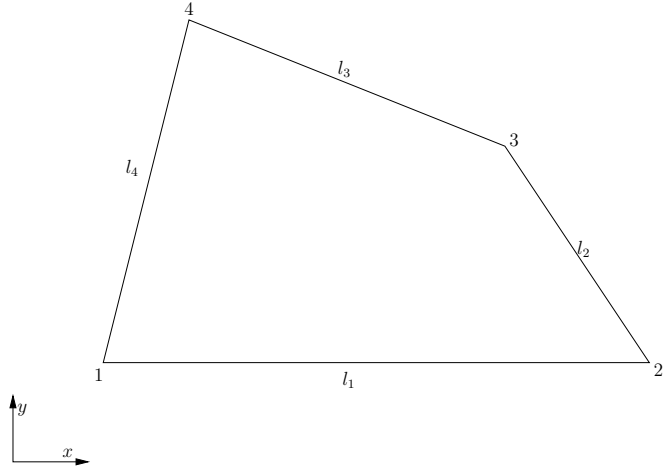}} \caption{ A sample
quadrilateral.} \label{fig:quad}
\end{figure}

\begin{subequations}
\begin{eqnarray}
l_{1}(x,y) = a_1x + b_1y + c_1 = 0 \\
l_{2}(x,y) = a_2x + b_2y + c_2 = 0 \\
l_{3}(x,y) = a_3x + b_3y + c_3 = 0 \\
l_{4}(x,y) = a_4x + b_4y + c_4 = 0
\end{eqnarray}
\end{subequations}

\noindent where $a_i,b_i$ and $c_i$ for $i=1,2,3,4$ are real constants. The wedge functions corresponding to each node, $w_i$, are defined so that they vary linearly along the edges adjacent to each node and vanish at the remaining nodes as~\cite{Wachspress1975}:

\begin{subequations}
\begin{eqnarray}
w_1(x,y) = \kappa_1~l_{2}(x,y)~l_{3}(x,y) \\
w_2(x,y) = \kappa_2~l_{3}(x,y)~l_{4}(x,y) \\
w_3(x,y) = \kappa_3~l_{4}(x,y)~l_{1}(x,y) \\
w_4(x,y) = \kappa_4~l_{1}(x,y)~l_{2}(x,y)
\end{eqnarray}
\end{subequations}

\noindent where the $\kappa_i$ are constants. In order that the Wachspress interpolants $N_i$, satisfy the partition of unity requirement, it is defined as:

\begin{equation}
N_i(x,y) = \frac{w_i}{\sum w_i}(x,y) %SPAB is the denominator nota function of x,y for Q4's? Is that what you mean by saying the functions are polynomials?
%Sundar: The denominator is evaluated at the node. So the denominator is a constant. For quadrailateral with opposite edges parallel, the external intersection point is at infinity. So the contribution is equal to 1. And the form, as given here is just polynomials.
\end{equation}

%\noindent In a more general form, the wedge functions can be expressed as:
%\begin{equation}
%w_i(x,y) = \Pi_{j \ne i,i+1}^{j=n}~l_{j}(x,y)
%\end{equation}

\begin{def1}
For a 4-sided polygon, the Wachspress rational basis interpolants degenerate to regular polynomials, i.e.  $\sum w_i(x,y)$ is a constant. %SPAB is this true?
\end{def1}

To illustrate the Wachspress interpolants for arbitrary quadrilaterals and to
answer Zhang {\it et al.,}~\cite{Zhang2008}'s queries on

\begin{figure} % Figure with averaged shape functions
\centering
\scalebox{1.0}{\input{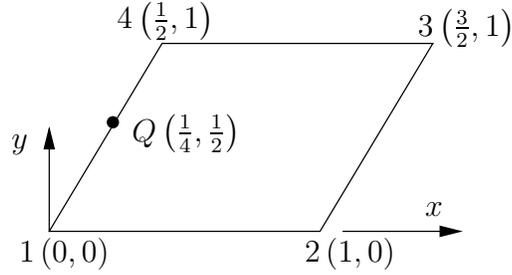}}
%\caption{A well-conditioned parallelogram}
\caption{ A parallelogram shaped element} \label{fig:parallelogram}
\end{figure}

\begin{itemize}
\item the existence of shape functions for arbitrary quadrilaterals;
\item positivity of the shape functions;
\item linearity of the shape functions,
\end{itemize}

\noindent we choose the same parallelogram element as used
in~\cite{Zhang2008} (cf. page 1292, Figure 2). This element is shown
in \fref{fig:parallelogram} for ease of reading.

\noindent The shape functions for the parallelogram shaped element (see \fref{fig:parallelogram})
based on Wachspress interpolation write:

\begin{subequations}
\begin{eqnarray}
N_1(x,y) = \left(y-1\right)\left(x-\frac{1}{2}y-1\right) \\
N_2(x,y) = \left(1-y\right)\left(x-\frac{1}{2}y\right) \\
N_3(x,y) = y\left(x-\frac{1}{2}y\right) \\
N_4(x,y) = -y\left(x-\frac{1}{2}y-1\right)
\end{eqnarray}
\end{subequations}

\noindent From the above equations, it can be easily verified that
the shape functions derived using the Wachspress approach satisfy
all the required properties (positivity, Kronecker delta property,
linear completeness) as shown in the book~\cite{Wachspress1975}.
What is more, the quality of these functions is independent of the
shape of the element, which make them ideal candidates for use in
the context of the SFEM where highly distorted meshes are of
interest given the absence of isoparametric mapping.

For the mid-point of side 1-4, i.e., the point $Q$ with coordinates $\left( \frac{1}{4},\frac{1}{2} \right)$ (see \fref{fig:parallelogram}), the shape function values are:

\begin{equation}
\renewcommand{\arraystretch}{1.5}
\capn(Q) = \left\{ \begin{array}{c} N_1(Q) \\ N_2(Q) \\ N_3(Q) \\ N_4(Q) \end{array} \right\} = \left\{ \begin{array}{c} \frac{1}{2} \\ 0 \\ 0 \\ \frac{1}{2} \end{array} \right\}
\end{equation}

\noindent as opposed to
\begin{equation}
\renewcommand{\arraystretch}{1.5}
\capn(Q) = \left\{ \begin{array}{c} N_1(Q) \\ N_2(Q) \\ N_3(Q) \\ N_4(Q) \end{array} \right\} = \left\{ \begin{array}{c} \frac{3}{8} \\ \frac{1}{8} \\ \frac{-1}{8} \\ \frac{5}{8} \end{array} \right\}
\end{equation}

\noindent as described in~\cite{Zhang2008} (cf. page. 1293).

\begin{def1} By a systematic approach using rational basis functions, the shape functions on arbitrary $n$-gons, in particular Q4's can be derived in the \emph{physical coordinate system}, without the need for isoparametric mapping, \emph{and} without the negative side effects of non-mapped Lagrange shape functions described in \cite{Zhang2008}. \end{def1}

\begin{def1} The form of the Wachspress shape functions depend on the element geometry and hence have to be recomputed for each element in the mesh, which is computationally expensive. However, an important advantage is that their quality does not deteriorate for highly distorted elements, including concave domains. Interested readers are referred to~\cite{Wachspress1975,Dasgupta2003,Sukumar2006} for advances in this direction. \end{def1}

%==================================================================
%=============== End of section on SFEM
%==================================================================

\section{NUMERICAL RESULTS}
\label{section:results}

\subsection{Patch Test}

It is straightforward to check that the linear patch test is
satisfied down to machine precision. The interested reader is
referred to the corresponding author to obtain a MATLAB code
including this and other test cases.

\subsection{Bending of a cantilever beam}
\begin{figure}
\centering
\includegraphics[angle=0,width=0.4\textwidth]{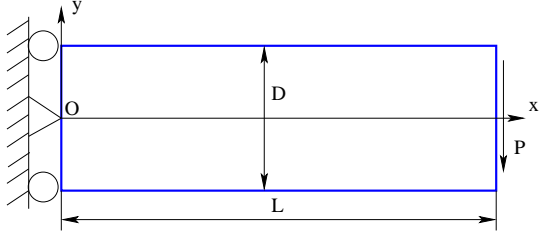}
\caption{Geometry and loads for a cantilever beam.}
\label{fig:beam}
\end{figure}

As a second example, bending of a thick cantilever subjected to a
parabolic load at the free end is examined as shown
in~\fref{fig:beam}. The geometry is: length $L$=8, height $D$=4 and
thickness $t$=1. The material properties are: Young's modulus
$E$=3$\times10^7$, and the parabolic shear force $P$=250. The exact
solution of this problem is available in~\cite{timo}.

%, giving the displacements in the $x$ and $y$ directions as

%\begin{subequations}
%\begin{eqnarray}
%\displaystyle
%u(x,y) &=& \frac{Py}{6\overline{E}I}\left[ (6L-3x)x+(2+\overline{\nu})(y^2-\frac{D^2}{4}) \right] \\
%\label{eqn:disp1} v(x,y) &=& -\frac{P}{6\overline{E}I}
%\left[3\overline{\nu}y^2(L-x)+(4+5\overline{\nu})\frac{D^2x}{4}+(3L-x)x^2
%\right]
%\end{eqnarray}
%\label{eqn:disp}
%\end{subequations}
%
%\noindent where $I$, the moment of inertia (second moment of area) is given
%by: \( I=\frac{D^3}{12}. \) and \(\overline{E}\) and
%\(\overline{\nu}\) in \Eref{eqn:disp} are given by
%
%\begin{equation}
%\renewcommand{\arraystretch}{1.5}
%    \overline{E} = \left\{ \begin{array}{ll}
%                    E & \mbox (\textup{plane\hspace{0.1cm}stress}), \\
%                    \frac{E}{1-\nu^2} & \mbox (\textup{plane\hspace{0.1cm}strain})
%                \end{array}
%                \right. \\ \hspace{0.25cm}
%    \overline{\nu} = \left\{ \begin{array}{ll}
%                        \nu & \mbox (\textup{plane\hspace{0.1cm}stress}), \\
%                        \frac{\nu}{1-\nu} & \mbox (\textup{plane\hspace{0.1cm}strain}).
%                    \end{array}
%                \right. \nonumber
%\end{equation}
%
%The stresses corresponding to the displacements in~\Eref{eqn:disp} are given by
%
%
%\begin{subequations}
%\begin{eqnarray}
%\displaystyle
%\sigma_{xx}(x,y)=\frac{P(L-x)y}{I} \\
%\sigma_{yy}(x,y) = 0 \\
%\sigma_{xy}(x,y) = -\frac{P}{2I} \left( \frac{D^2}{4}-y^2 \right)
%\end{eqnarray}
%\end{subequations}

In this problem, two types of mesh are considered: one is uniform and regular and the other is an
irregular mesh, with the coordinates of interior nodes given by:

\begin{subequations}
\begin{eqnarray}
x^{\prime} = x + \left( 2r_c-1 \right) \alpha_{ir} \Delta x \\
y^{\prime} = y + \left( 2r_c-1 \right) \alpha_{ir} \Delta y \\
\end{eqnarray}
\end{subequations}

\noindent where $r_c$ is a random number between $0$ and $1.0$, $\alpha_{ir} \in \left[0,0.5\right]$ is
an irregularity factor that controls the shapes of the elements and $\Delta x$ and $\Delta y$ are the
initial regular element sizes in the $x-$ and $y-$ directions respectively.

\begin{figure}
\centering
\subfigure[]{\includegraphics[angle=0,width=0.7\textwidth]{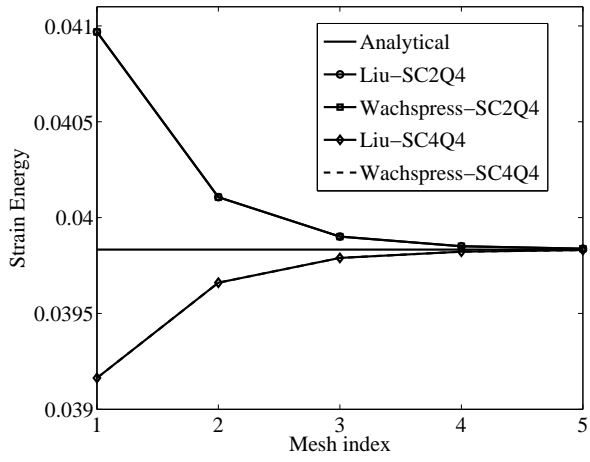}}
\subfigure[]{\includegraphics[angle=0,width=0.7\textwidth]{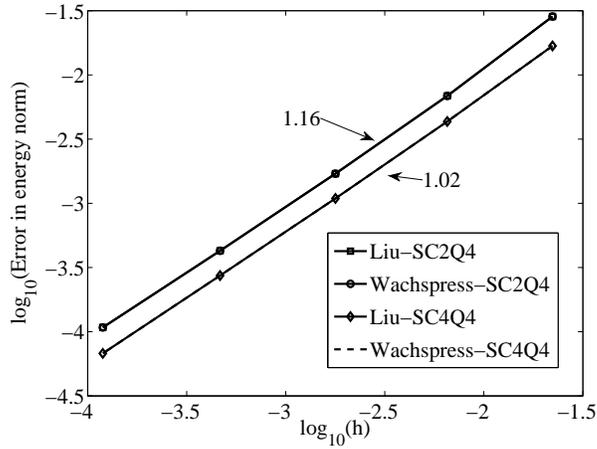}}
\caption{The convergence of the numerical energy to the exact energy
and convergence rate in the energy norm with regular meshes for the
cantilever beam problem: (a) strain energy and (b) the convergence
rate.} \label{fig:seregular}
\end{figure}

\begin{figure}
\centering
\subfigure[]{\includegraphics[angle=0,width=0.68\textwidth]{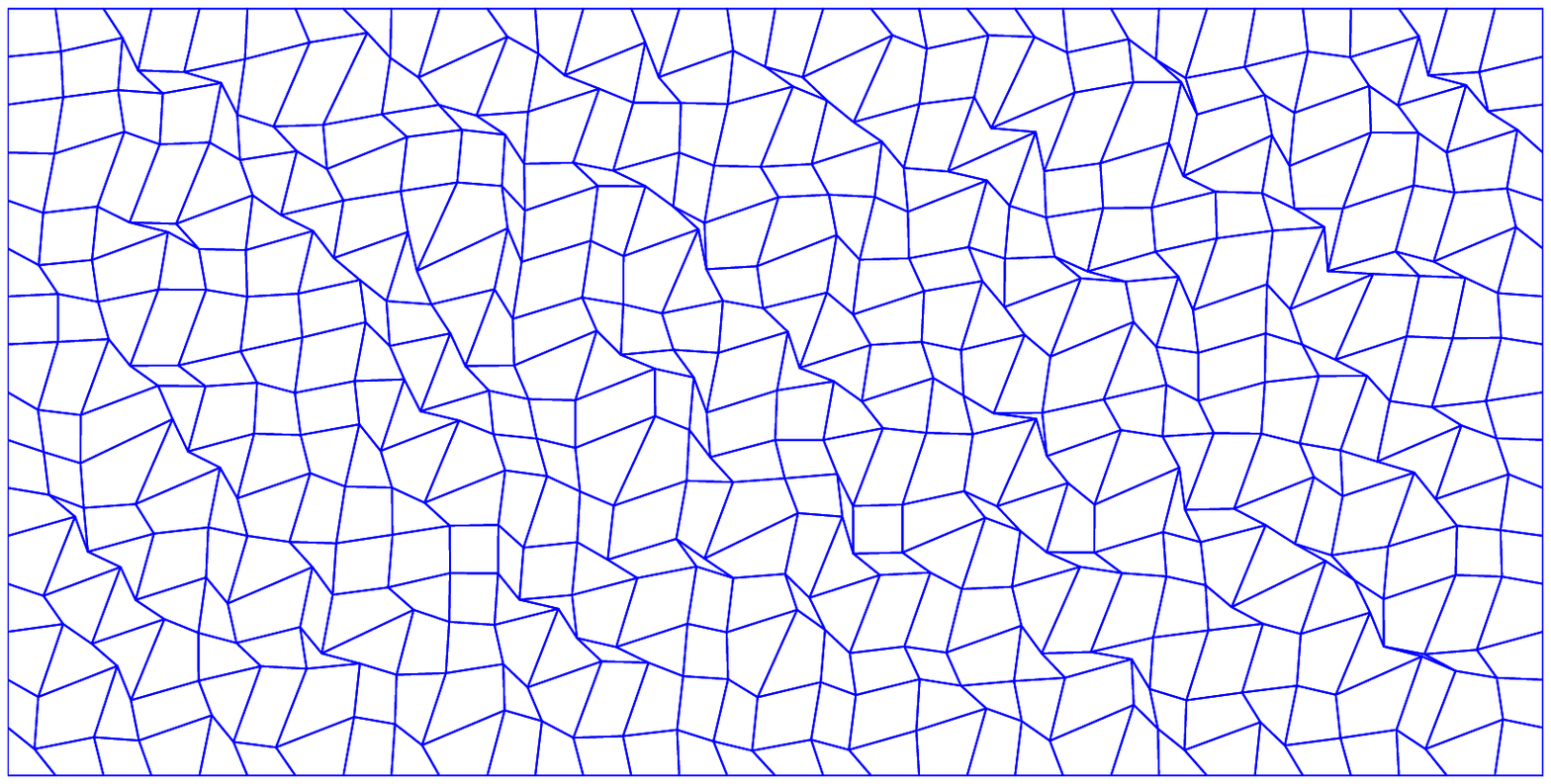}}
\subfigure[]{\includegraphics[angle=0,width=0.68\textwidth]{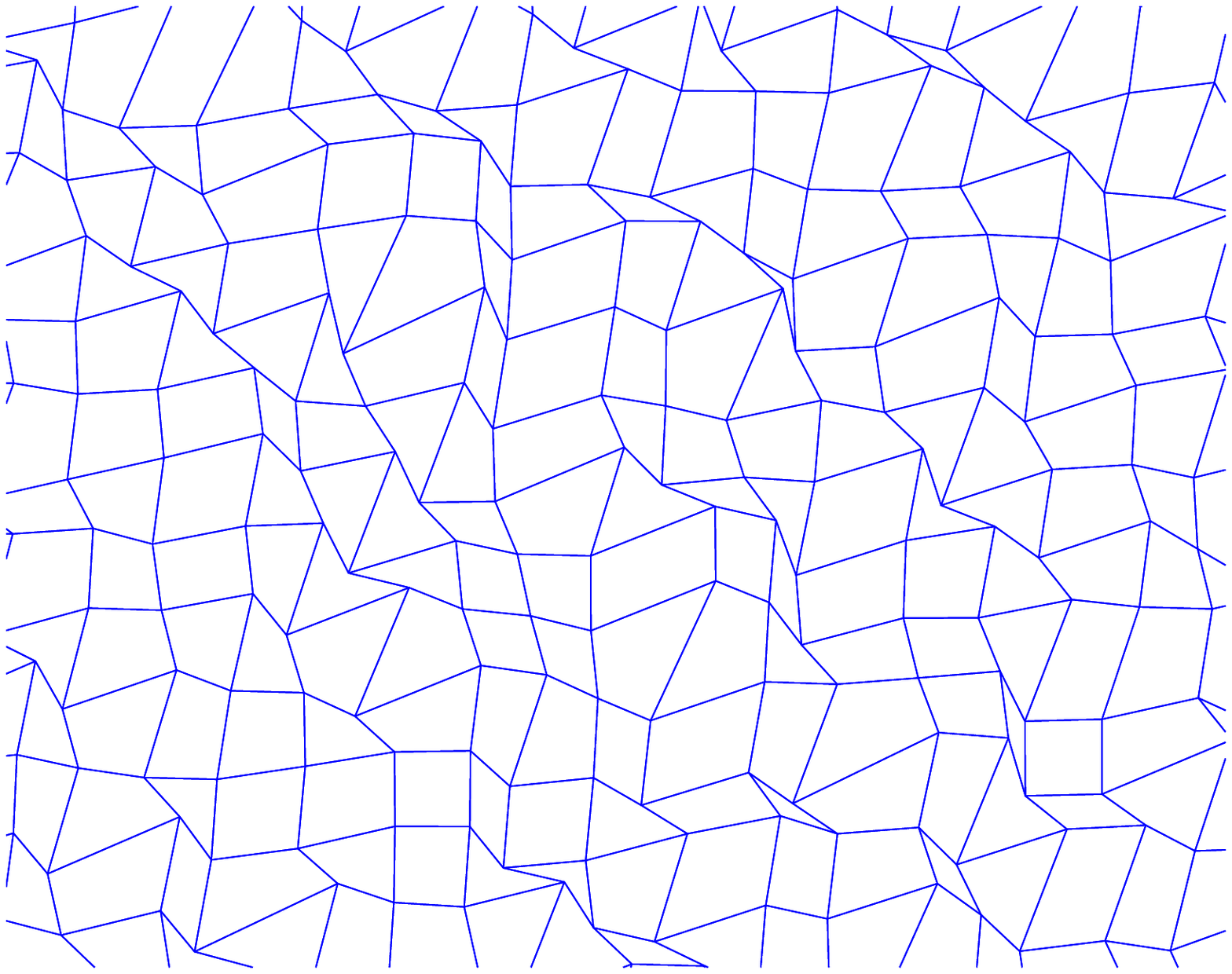}}
\caption{Cantilever beam:(a) irregular mesh with extremely distorted
elements and (b) zoomed view of the mesh.} \label{fig:meshirr}
\end{figure}

\begin{figure}
\centering
\includegraphics[angle=0,width=0.7\textwidth]{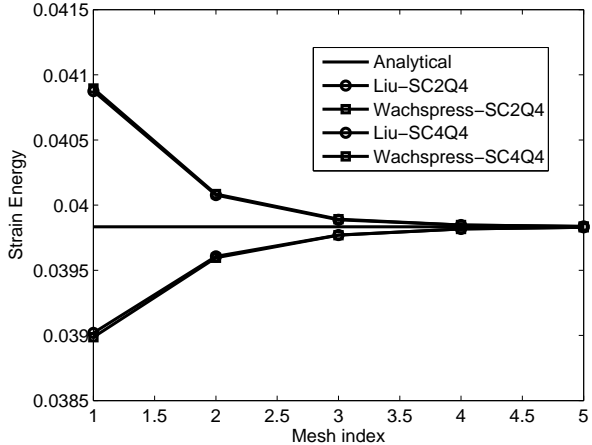}
\caption{The convergence of the numerical energy to the exact energy
with irregular meshes
$(\alpha_{ir}=0.2)$ for the cantilever beam problem.} \label{fig:seirregular}
\end{figure}

\begin{figure}
\centering
\subfigure[]{\includegraphics[angle=0,width=0.7\textwidth]{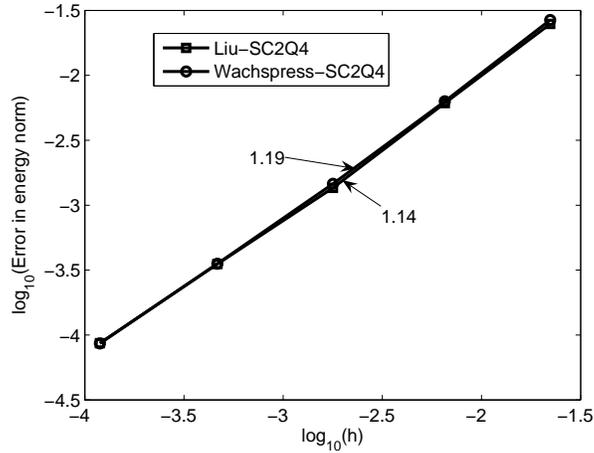}}
\subfigure[]{\includegraphics[angle=0,width=0.7\textwidth]{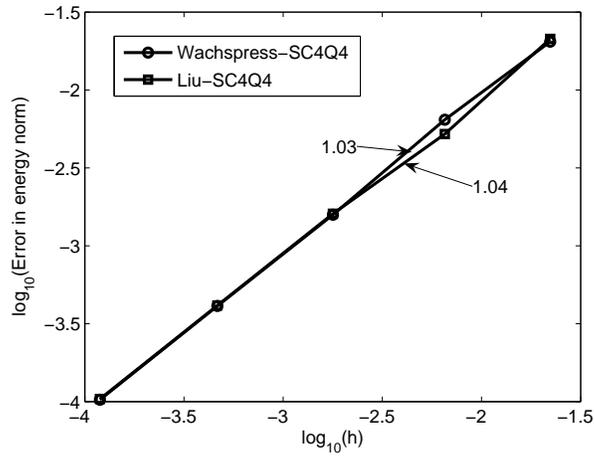}}
\caption{The convergence rate in the energy norm with irregular meshes
$(\alpha_{ir}=0.2)$ for the cantilever beam problem: (a) SC2Q4 and (b) SC4Q4.} \label{fig:seirregular1}
\end{figure}

\begin{figure}
\centering
\subfigure[]{\includegraphics[angle=0,width=0.7\textwidth]{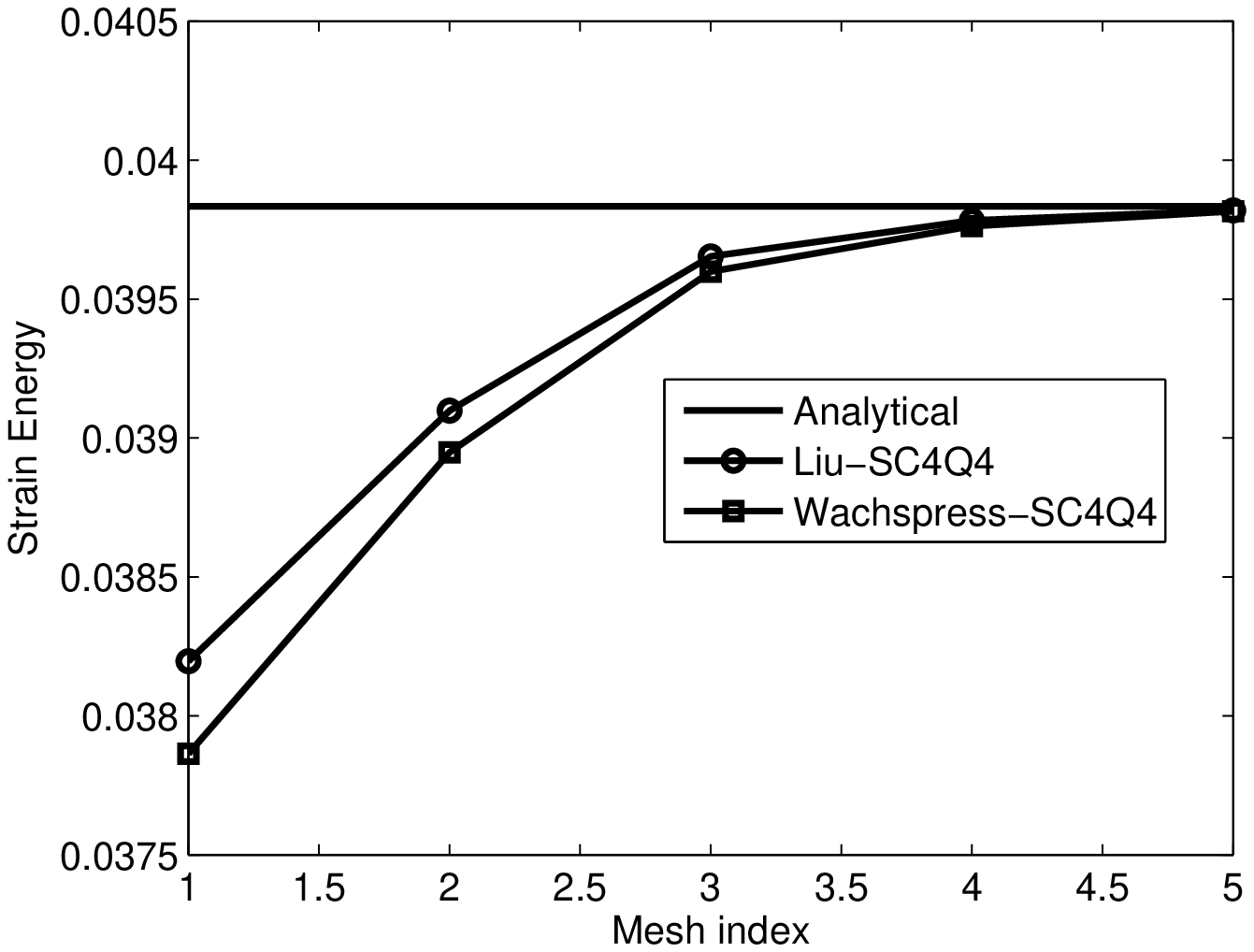}}
\subfigure[]{\includegraphics[angle=0,width=0.7\textwidth]{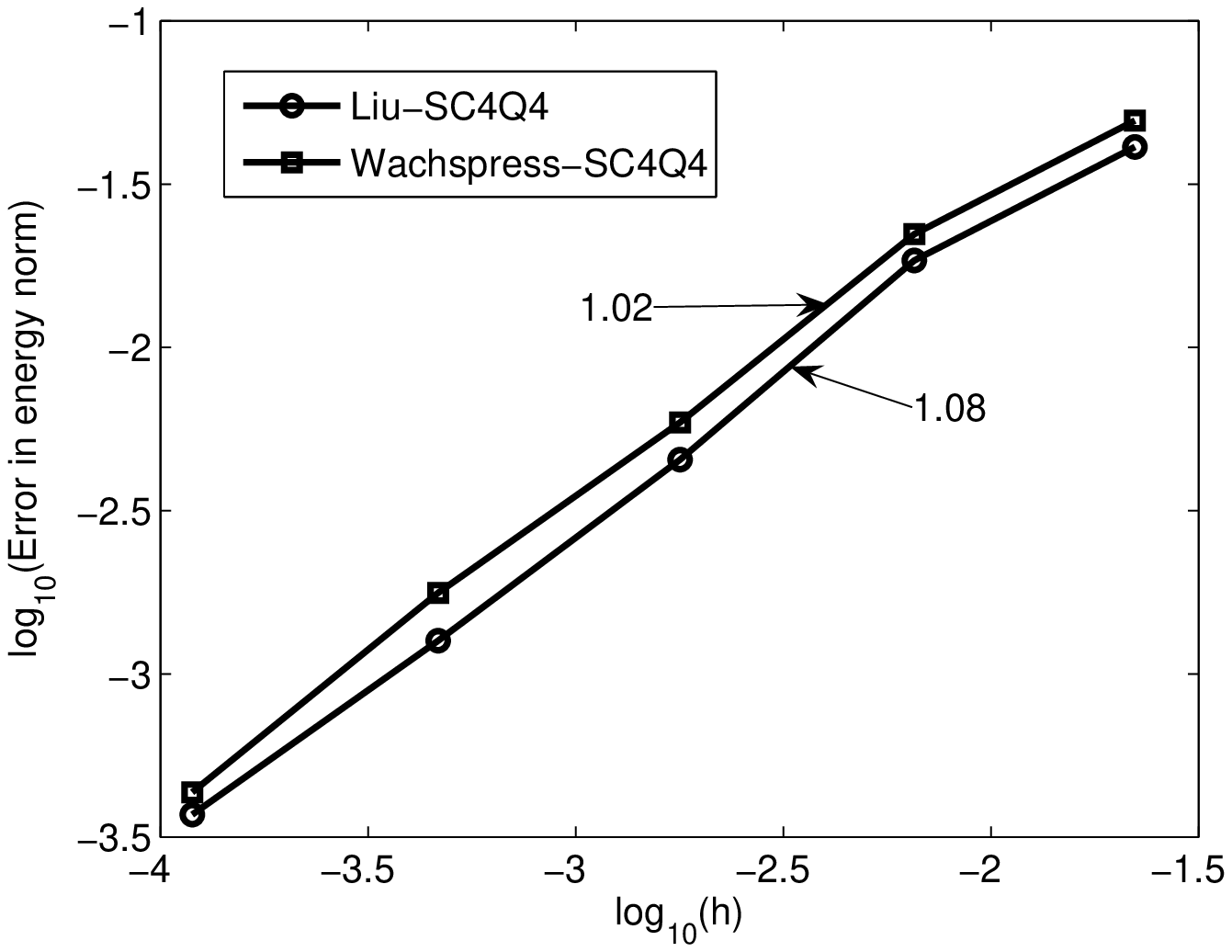}}
\caption{The convergence of the numerical energy to the exact energy
and convergence rate in the energy norm with irregular meshes
$(\alpha_{ir}=0.5)$ for the cantilever beam problem: (a) strain
energy and (b) the convergence rate.} \label{fig:seirregular2}
\end{figure}

The domain is discretized with quadrilateral elements and two
quadrilateral subcells for each element are used for the current
study, denoted by $SCkQ4$\footnote{SCkQ4 implies k quadrilateral
subcells for 4 noded quadrilateral elements}. \fref{fig:meshirr}
shows a discretization of the cantilever beam with quadrilateral
elements using an irregular mesh. A value of $\alpha_{ir}=0.5$ is
used to generate the irregular mesh. Under plane stress conditions
and for a Poisson ratio $\nu$=0.3, the exact strain energy is
0.0398333. \fref{fig:seregular} illustrates the convergence of the
strain energy and the rate of convergence in the energy norm of the
elements built using the Wachspress interpolants compared with those
of the SFEM results (\cite{hbb0}) for regular quadrilateral meshes.
The total strain energy and the rate of convergence for different irregular
meshes $(\alpha_{ir}=0.2 ~\textup{and} ~\alpha_{ir}=0.5)$ are shown in \frefs{fig:seirregular}, (\ref{fig:seirregular1}) and
(\ref{fig:seirregular2}).

%----- irregular mesh results

It is seen that for regular meshes, the results are identical and for irregular meshes, there is a small difference, but the rate of convergence is not affected significantly.
It is seen that with increase in mesh index\footnote{$\textup{Mesh index} = \frac{\textup{number of
elements in the $x$-direction}}{\textup{length of the cantilever
beam}}$}, both  methods converge to the exact solution. With the use of Wachspress interpolants, the question of positivity of the shape functions and existense of the shape functions for arbitrary quadrilaterals is resolved yet preserving the true essence of the SFEM.

%==================================================================
%=============== End of section on Numerical Results
%==================================================================
\section{CONCLUSIONS}

This letter showed that it is possible to both retain the highly desirable features of the smoothed FEM (SFEM), and its true essence, i.e. strain smoothing and boundary integration, without sacrificing a rigorous approximation where the shape functions are known explicitly at any point of the smoothed finite element.

The proposed interpolation scheme, known as Wachspress interpolation, provides a suitable means to suppress the problems of definition, positivity and linearity associated with non-mapped Lagrange shape functions. It also provides a general framework to define approximation over arbitrary (possibly non-convex) polytopes including curved edges or surfaces.

The Wachspress approximation  was conclusively tested within the context of the smoothed FEM, showing that it yields a method which passes the patch test and provides, even for very high element distortion, accurate and optimally convergent results in linear elasticity.

The results also showed that utilizing the proposed Wachspress shape functions to build the approximation in the smoothed FEM leads to no significant difference in results compared to the widely used technique of ``averaged shape functions'' preconized in the original paper \cite{Liu2007} and used in all the smoothed FEM work published to date \cite{Hung2009,hbb_shell,Nguyen2007,hbb0,hbb_plate,liu:sfem2,hbb_shell}.

This finding reinforces the claim made in \cite{Liu2009} that the approximation need not be known explicitly to build an optimal SFEM and dispels, also, the misconception of \cite{Zhang2008}, whose detraction of the smoothed FEM rests on the misunderstanding that this method requires the use of non-mapped Lagrange shape functions.

It will be most interesting to investigate the behavior of the smoothed FEM, especially the choice of approximation, in the context of higher-order methods as well as non-polynomial enrichments \cite{Bordas2008}.

{\bf Acknowledgement}
The authors would like to thank Professor Ted Belytschko for his expert comments, which were very useful in improving the quality, clarity and soundness of the paper. The support of the Overseas Research Students Awards Scheme and of the Faculty of Engineering, University of Glasgow is gratefully acknowledged.

\bibliographystyle{wileyj}
\bibliography{sfem_reference,mesh_free_reference,PolyFEM,fem_reference}

\begin{thebibliography}{10}
\providecommand{\url}[1]{\texttt{#1}}
\providecommand{\urlprefix}{URL }
\expandafter\ifx\csname urlstyle\endcsname\relax
  \providecommand{\doi}[1]{doi:\discretionary{}{}{}#1}\else
  \providecommand{\doi}{doi:\discretionary{}{}{}\begingroup
  \urlstyle{rm}\Url}\fi

\bibitem{Zhang2008}
Zhang HH, Liu SJ, Li LX. On the smoothed finite element method.
  \emph{International Journal of Numerical Methods in Engineering}  2008;
  \textbf{76}(8):1285--1295, \doi{10.1002/nme.2460}.

\bibitem{Liu2009}
Liu G, Nguyen-Thoi T, Nguyen-Xuan H, Dai K, Lam K. On the essence and the
  evaluation of the shape functions for the smoothed finite element method
  (sfem). \emph{International Journal of Numerical Methods in Engineering}
  2009; \doi{10.1002/nme.2587}.

\bibitem{Liu2007}
Liu GR, Dai KY, Nguyen TT. A smoothed finite element for mechanics problems.
  \emph{Computational Mechanics}  May 2007; \textbf{39}(6):859--877,
  \doi{10.1007/s00466-006-0075-4}.

\bibitem{Hung2009}
Hung NX, Bordas S, Hung N. Addressing volumetric locking and instabilities by
  selective integration in smoothed finite element. \emph{Communications in
  Numerical Methods in Engineering}  2009; \textbf{25}(1):19--34,
  \doi{10.1002/cnm.1098}.

\bibitem{liu:sfem2}
Liu GR, Nguyen TT, Dai KY, Lam KY. Theoretical aspects of the smoothed finite
  element method ({SFEM}). \emph{Int. J. Numer. Meth. Engng.}  2007;
  \textbf{71}(8):902--930.

\bibitem{hbb_shell}
Nguyen NT, Rabczuk T, Nguyen-Xuan H, Bordas S. A smoothed finite element method
  for shell analysis. \emph{Computer Methods in Applied Mechanics and
  Engineering}  2008; \textbf{198}(2):165--177,
  \doi{10.1016/j.cma.2008.05.029}.

\bibitem{Nguyen2007}
Nguyen XH, Bordas S, Nguyen-Dang H. Deduction of pure displacement and
  equilibrium elements by smoothed integration technique. \emph{Computer
  Methods in Applied Mechanics and Engineering}  2007; .

\bibitem{hbb0}
Nguyen-Xuan H, Bordas S, Nguyen-Dang H. Smooth finite element methods:
  Convergence, accuracy and properties. \emph{Int. J. Numer. Meth. Engng.}
  2008; \textbf{74}(2):175--208, \doi{10.1002/nme.2146}.

\bibitem{hbb_plate}
Nguyen-Xuan H, Rabczuk T, Bordas S, Debongnie JF. A smoothed finite element
  method for plate analysis. \emph{Computer Methods in Applied Mechanics and
  Engineering}  2008; \textbf{197}(13-16):1184--1203,
  \doi{10.1016/j.cma.2007.10.008}.

\bibitem{Wachspress1975}
Wachspress EL. \emph{A rational basis for function approximation}. Academic
  Press, Inc. New York., 1975.

\bibitem{Chen2000}
Chen JS, Wang HP. Some recent improvements in meshfree methods for
  incompressible finite elasticity boundary value problems with contact.
  \emph{Computational Mechanics}  2000; \textbf{25}:137--156.

\bibitem{Bordas2008}
Bordas SP, Rabczuk T, Hung NX, Nguyen VP, Natarajan S, Bog T, Quan DM, Hiep NV.
  Strain smoothing in fem and xfem. \emph{Computers and Structures}  2008;
  \doi{10.1016/j.compstruc.2008.07.006}.

\bibitem{Dasgupta2003}
Dasgupta G. Interpolants with convex polygons: Wachspress's shape functions.
  \emph{Journal of Aerospace Engineering}  January 2003; \textbf{16}(1):1--8.

\bibitem{Sukumar2006}
Sukumar N, Malsch EA. Recent advances in the construction of polygonal finite
  element interpolants. \emph{Archives of Computational Methods in Engineering}
   2006; \textbf{13}(1):129--163.

\bibitem{timo}
Timoshenko SP, Goodier JN. \emph{Theory of elasticity}. Mc-Graw-Hill, New York,
  1970.

\end{thebibliography}

\end{document}